# On classifying decision making units in DEA: A unified dominance-based model

Mahmood Mehdiloozad

*Department of Mathematics, College of Sciences, Shiraz University 71454, Shiraz, Iran*

Mohammad Bagher Ahmadi [*]

*Department of Mathematics, College of Sciences, Shiraz University 71454, Shiraz, Iran*

Biresh K. Sahoo

*Xavier Institute of Management, Xavier University, Bhubaneswar 751 013, India*

[*] **Corresponding author**: M. B. Ahmadi

Associate Professor of Mathematics

Department of Mathematics

College of Sciences

Shiraz University

Golestan Street | Adabiat Crossroad | Shiraz 71454 | Iran

E-mail: mbahmadi@shirazu.ac.ir

Tel.: +98-917-717-0280

# On classifying decision making units in DEA: A unified dominance-based model

Abstract

In data envelopment analysis (DEA), the concept of efficiency is examined in either Farrell (DEA) or Pareto senses. In either of these senses, the efficiency status of a decision making unit (DMU) is classified as either *weak* or *strong*. It is well established that the strong DEA efficiency is both necessary and sufficient for achieving the Pareto efficiency. For the weak Pareto efficiency, however, the weak DEA efficiency is only sufficient, but not necessary in general. Therefore, a DEA-inefficient DMU can be either weakly Pareto efficient or Pareto inefficient. Motivated by this fact, we propose a new classification of DMUs in terms of both DEA and Pareto efficiencies. To make this classification, we first demonstrate that the Farrell efficiency is based on the notion of FGL dominance. Based on the concept of dominance, we then propose and substantiate an alternative single-stage method. Our method is computationally efficient since (i) it involves solving a unique single-stage model for each DMU, and (ii) it accomplishes the classification of DMUs in both input and output orientations simultaneously. Finally, we present a numerical example to illustrate our proposed method.

## 1 Introduction

Data envelopment analysis (DEA), developed by Charnes et al. (1978), deals with measuring the relative efficiency of a homogeneous group of observed decision making units (DMUs). Such DMUs are transforming, in varying quantities, the same set of inputs into the same set of outputs. Based on observed input–output data and a set of postulates, DEA defines a reference technology relative to which the efficiency of each individual observation is estimated.

In the DEA literature, the concept of efficiency has been considered and analyzed in two senses. The first one, on which the pioneering CCR model of Charnes et al. (1978) and the BCC model of Banker et al. (1984) are based, is due to Debreu (1951) and Farrell (1957)[1]. Under variable returns to scale (VRS) framework, the Debreu–Farrell efficiency can be measured using the BCC model. With an input orientation, a DMU is *Debreu–Farrell efficient* or *weakly BCC efficient* if and only if the pro rata improvements to inputs are not possible; otherwise, it is called *BCC inefficient*. Moreover, a weakly BCC-efficient DMU is *strongly BCC efficient* when no non-pro rata improvements to inputs or outputs are possible without compromising technological feasibility. Note that a strongly BCC-efficient DMU is either *extreme* or *non-extreme* depending on whether it is represented as a convex combination of the remaining DMUs. The second one, which is attributable to Pareto (1909) and Koopmans (1951), is based on the concept of *dominance*. A DMU is *(weakly) strongly Pareto–Koopmans efficient* if and only if no other technologically feasible activity (strongly) dominates it; otherwise, it is *Pareto–Koopmans inefficient*. For brevity, we henceforth refer to Debreu–Farrell and Pareto–Koopmans efficiencies as Farrell and Pareto efficiencies, respectively.

It is well known that the strong BCC and Pareto efficiencies are equivalent (Cooper et al. 2007). As regards the relationships between the weakly BCC-efficient, weakly Pareto-efficient and boundary activities, Krivonozhko et al. (2005) established some new interesting findings. First, the input- (output-) oriented weak BCC efficiency implies the weak Pareto efficiency, but the reverse is not generally true. Second, while the weak Pareto efficiency is both necessary and sufficient for a DMU to be on the boundary, the weak BCC efficiency is only sufficient but not necessary. As a consequence, it is not guaranteed that a BCC-inefficient unit is in the interior of the technology set, and it may, hence, be a boundary activity.

An important question arises now as to how to make a classification of DMUs in terms of the BCC and Pareto efficiencies. To our knowledge, a few research studies have been devoted to the BCC-efficiency based classification (see, e.g., Charnes et al. (1986, 1991), Thrall (1996) and Dulá and Hickman (1997), among others). These methods partition the observed DMUs into four groups: extreme efficient, non-extreme efficient, weakly BCC efficient and BCC inefficient. Though the existing classification methods are interesting, as we will be shortly demonstrating, they fail to identify whether a BCC-inefficient unit is either an interior or a boundary activity. On the other hand, to the best of our knowledge, barring one study by Jahanshahloo et al. (2005), no other research works have so far been conducted on classifying DMUs based on the Pareto efficiency.

---

[1] Throughout our paper, we define the concept of Farrell efficiency with respect to *strongly* disposable technology and refer to it as the *Farrell* measure of efficiency. However, as per Färe et al. (1985), it is called as the *weak* measure of efficiency, as they define the Farrell measure of efficiency with respect to weakly disposable technology.

In this contribution, we propose a new classification method that simultaneously classifies DMUs in terms of both BCC and Pareto efficiencies, but without resorting to their efficiency scores. To pursue this objective, we first make a reclassification of all the observed DMUs by dividing the BCC-inefficient units into the interior and boundary activities. Although making this classification via the super-efficiency models of Andersen and Petersen (1993) and Jahanshahloo et al. (2005) is feasible, but it requires solving a large number of problems. We therefore propose, based on the concept of dominance, an alternative single-stage model that enables us to simultaneously make our proposed classification in both input and output orientations.

As the first step of developing our model, we establish a relationship between the Farrell efficiency and the dominance concept. Specifically, we demonstrate that the Farrell efficiency is based on the notion of *FGL dominance*[2]. Analogous to the procedure followed in the super-efficiency measurement, we exclude the DMU under evaluation from the reference set and then, construct the technology set with the remaining DMUs. We define the set of activities in this technology set that all dominate the DMU under evaluation. We then prove that the problem under consideration reduces, to checking whether the defined set is empty, and if not, to identifying an activity in this set such that its corresponding input–output slack vector has the maximum number of positive components. Based on this finding, we propose our classification model.

Our proposed method has three advantages. First, it involves solving of a unique single linear programming problem. Second, our model classifies DMUs in terms of input and output orientations simultaneously. With regard to these two advantages, our approach can be considered computationally more efficient than the excising ones. Third, by identifying all the inputs and outputs that need to be improved, our proposed model enables a decision maker to deal with his/her preference as to which specific input to be reduced or which specific output to be increased by a weakly Pareto-efficient DMU in order to improve its overall performance. We note that since our model contains several upper-bounded variables, the computational efficiency of our method can be enhanced by using the simplex algorithm[3] adopted for solving the LP problems with upper-bounded variables.

---

[2] Thought Färe et al. (1985, pp. 28 & 46) employ this notion to characterize the weakly efficient subsets of the input and output correspondences and the graph of the technology; they do not call it as a dominance concept.

The remainder of this paper unfolds as follows. Section 2 discusses the background of our research. Section 3 presents the main contribution of our study, where a unified classification model is proposed. Section 4 illustrates the use of our proposed method with two numerical examples. Section 5 concludes with some remarks.

## 2 Background of the research

### *2.1 Technology set*

As far as notations are concerned, let $\mathbb{R}_+^d$ be the non-negative Euclidean *d*-orthant. We denote vectors and matrices in bold letters, vectors in lower case and matrices in upper case. All vectors are assumed to be column vectors. We denote by a subscript $T$ the transpose of vectors and matrices. We also use $\mathbf{0}_d$ and $\mathbf{1}_d$ to show *d*-dimensional vectors with the values of 0 and 1 in every entry, respectively. For $\mathbf{u}, \mathbf{v} \in \mathbb{R}^d$, we use the following standard notations that are taken from Färe et al. (1985):

$$\begin{aligned} \mathbf{u} \geqq \mathbf{v} &\Leftrightarrow u_j \geq v_j \text{ for all } j = 1,...,d; \\ \mathbf{u} \geq \mathbf{v} &\Leftrightarrow u_j \geq v_j \text{ for all } j = 1,...,d \text{ and } \mathbf{u} \neq \mathbf{v}; \\ \mathbf{u} > \mathbf{v} &\Leftrightarrow u_j > v_j \text{ for all } j = 1,...,d. \end{aligned} \quad (1)$$

If $\mathbf{u}, \mathbf{v} \in \mathbb{R}_+^d$, we write $\mathbf{u} >^+ \mathbf{v}$ to denote that $u_j > v_j$ or $u_j = v_j = 0$ for all $j = 1,...,d$.

We consider $O$ as a set of *n* observed DMUs, where each uses *m* inputs to produce *s* outputs. For any $j \in \{1,...,n\}$, we denote the *j*th DMU by DMU$_j$ and its input and output vectors by $\mathbf{x}_j \in \mathbb{R}_+^m$ and $\mathbf{y}_j \in \mathbb{R}_+^s$, respectively. The *i*th ($i = 1,...,m$) input and the *r*th ($r = 1,...s$) output of DMU$_j$ are symbolized by $x_{ij}$ and $y_{rj}$, respectively.

The *technology set* (*T*) is defined as the set of all feasible input–output combinations $(\mathbf{x}, \mathbf{y})$, i.e.,

$$T = \left\{ (\mathbf{x}, \mathbf{y}) \in \mathbb{R}_+^m \times \mathbb{R}_+^s \mid \mathbf{x} \text{ can produce } \mathbf{y} \right\}. \quad (2)$$

---

[3] The simplex algorithm for bounded variables was published by Dantzig (1955) and was independently developed by Charnes and Lemke (1954). The use of this algorithm is much more efficient than the ordinary simplex algorithm for solving the LP problem with upper-bounded variables (Winston 2003).

Following Banker et al. (1984) and Andersen and Petersen (1993), we define $T^{DEA}$ and $T_o^{DEA}$ as the technology sets generated, respectively, by the observed DMUs and by the observed DMUs excluding $DMU_o$. Then, under variable returns to scale (VRS) framework, the nonparametric DEA-based representations $T^{DEA}$ and $T_o^{DEA}$ are, respectively, set up as

$$T^{DEA} = \left\{ (\mathbf{x},\mathbf{y}) \in \mathbb{R}_+^m \times \mathbb{R}_+^s \,\middle|\, \mathbf{x} \geq \mathbf{X}\boldsymbol{\lambda},\ \mathbf{Y}\boldsymbol{\lambda} \geq \mathbf{y},\ \mathbf{1}_n^T \boldsymbol{\lambda} = 1,\ \boldsymbol{\lambda} \geq \mathbf{0}_n \right\}, \tag{3}$$

$$T_o^{DEA} = \left\{ (\mathbf{x},\mathbf{y}) \in \mathbb{R}_+^m \times \mathbb{R}_+^s \,\middle|\, \mathbf{x} \geq \mathbf{X}_o \boldsymbol{\mu},\ \mathbf{Y}_o \boldsymbol{\mu} \geq \mathbf{y},\ \mathbf{1}_{n-1}^T \boldsymbol{\mu} = 1,\ \boldsymbol{\mu} \geq \mathbf{0}_{n-1} \right\}. \tag{4}$$

The subscript $o \in \{1,...,n\}$ is the index of the DMU under evaluation. Here, $\mathbf{X} = [\mathbf{x}_1 \ ... \ \mathbf{x}_n]$ and $\mathbf{Y} = [\mathbf{y}_1 \ ... \ \mathbf{y}_n]$ represent the input and output matrices, respectively; further, $\mathbf{X}_o$ and $\mathbf{Y}_o$ are, respectively, obtained from $\mathbf{X}$ and $\mathbf{Y}$ by removing the columns $\mathbf{x}_o$ and $\mathbf{y}_o$.

## *2.2 Efficiency in the senses of Farrell and Pareto*

With reference to $T^{DEA}$, the input-oriented BCC model of Banker et al. (1984) can be formulated as

$$\begin{aligned}
&\min\quad \theta - \varepsilon \left( \mathbf{1}_m^T \mathbf{s}^- + \mathbf{1}_s^T \mathbf{s}^+ \right) \\
&\text{subject to} \\
&-\mathbf{X}\boldsymbol{\lambda} - \mathbf{s}^- \geq -\theta \mathbf{x}_o, \\
&\mathbf{Y}\boldsymbol{\lambda} - \mathbf{s}^+ \geq \mathbf{y}_o, \\
&\mathbf{1}_n^T \boldsymbol{\lambda} = 1, \\
&\boldsymbol{\lambda} \geq \mathbf{0}_n,\ \mathbf{s}^- \geq \mathbf{0}_m,\ \mathbf{s}^+ \geq \mathbf{0}_s,
\end{aligned} \tag{5}$$

where $\varepsilon > 0$ is a positive non-Archimedean infinitesimal; and, $\mathbf{s}^-$ and $\mathbf{s}^+$ are the input and output slack vectors, respectively.

Note that the use of $\varepsilon$ in the objective function of model (5) enables us to detect the presence of non-zero slacks. When the value of $\varepsilon$ is specified, model (5) can be solved in a single stage. It has, however, been pointed out that the single-stage approach may result in computational inaccuracies that lead to results due to the choice of $\varepsilon$ (Ali and Seiford 1993, Chang and Guh 1991). Hence, the standard approach for avoiding the use of $\varepsilon$ in practical applications is to apply the following two-stage procedure. Stage 1: minimize $\theta$ by ignoring the slacks. Stage 2: replace $\theta$ by its optimal value ($\theta^*$) in (5), and maximize the sum of the slacks.

**Definition 2.1** Let $\left(\theta^*, \boldsymbol{\lambda}^*, \mathbf{s}^{-*}, \mathbf{s}^{+*}\right)$ be an optimal solution to model (5). Then, DMU$_o$ is called

- *BCC inefficient* if and only if $\theta^* < 1$.
- *weakly BCC efficient* if and only if $\theta^* = 1$.
- *strongly BCC efficient* if and only if $\theta^* = 1$, $\mathbf{s}^{-*} = \mathbf{0}_m$ and $\mathbf{s}^{+*} = \mathbf{0}_s$.

With respect to $T^{DEA}$, we define the concept of dominance as follows.

**Definition 2.2** Let $(\mathbf{x}', \mathbf{y}')$ and $(\mathbf{x}, \mathbf{y})$ be two input–output combinations in $T^{DEA}$. Then,

- $(\mathbf{x}', \mathbf{y}')$ *(weakly) dominates* $(\mathbf{x}, \mathbf{y})$ in Pareto sense, i.e., $(\mathbf{x}', \mathbf{y}') \succeq_P (\mathbf{x}, \mathbf{y})$, if and only if $\begin{pmatrix} -\mathbf{x}' \\ \mathbf{y}' \end{pmatrix} \geq \begin{pmatrix} -\mathbf{x} \\ \mathbf{y} \end{pmatrix}$.

- $(\mathbf{x}', \mathbf{y}')$ *strongly dominates* $(\mathbf{x}, \mathbf{y})$ in Pareto sense, i.e., $(\mathbf{x}', \mathbf{y}') \succ_P (\mathbf{x}, \mathbf{y})$, if and only if $\begin{pmatrix} -\mathbf{x}' \\ \mathbf{y}' \end{pmatrix} > \begin{pmatrix} -\mathbf{x} \\ \mathbf{y} \end{pmatrix}$.

The above definition of dominance leads to the following definition of efficiency in the sense of Pareto.

**Definition 2.3** DMU$_o$ is called

- *Pareto inefficient* if and only if $(\mathbf{x}, \mathbf{y}) \succ_P (\mathbf{x}_o, \mathbf{y}_o)$ holds for some $(\mathbf{x}, \mathbf{y}) \in T^{DEA}$.
- *weakly Pareto efficient* if and only if $(\mathbf{x}, \mathbf{y}) \succ_P (\mathbf{x}_o, \mathbf{y}_o)$ holds for no $(\mathbf{x}, \mathbf{y}) \in T^{DEA}$.
- *strongly Pareto efficient* if and only if $(\mathbf{x}, \mathbf{y}) \succeq_P (\mathbf{x}_o, \mathbf{y}_o)$ holds for no $(\mathbf{x}, \mathbf{y}) \in T^{DEA}$.

The following proposition, taken from Cooper et al. (2007), demonstrates that the strong BCC and Pareto efficiencies are equivalent.

**Proposition 2.1** DMU$_o$ is strongly BCC efficient if and only if it is strongly Pareto efficient.

As regards the relationships between the weakly DEA-efficient, weakly Pareto-efficient and boundary activities, Krivonozhko et al. (2005) established some interesting findings in their impressive study that had not been considered theretofore in the DEA literature.

**Proposition 2.2** (Krivonozhko et al. 2005) $DMU_o$ is weakly Pareto efficient if and only if it belongs to $\partial T^{DEA}$ (the boundary of $T^{DEA}$).

The above proposition states that $\partial T^{DEA}$ is characterized by the weakly Pareto-efficient DMUs. As an immediate result, $DMU_o$ is Pareto inefficient if and only if it lies in $\text{int}(T^{DEA})$ (the interior of $T^{DEA}$).

**Proposition 2.3** (Krivonozhko et al. 2005) If $DMU_o$ is weakly BCC efficient in either input orientation or output orientation or both, then it belongs to $\partial T^{DEA}$.

Note that the reverse of Proposition 2.3 is not true in general. In order to illustrate this, let us consider unit E in Fig. 2 which is a boundary activity. From Proposition 2.2, this unit is weakly Pareto efficient. However, it is not BCC efficient since its inputs can be decreased radially. Thus, while the weak Pareto efficiency is both necessary and sufficient for a DMU to be a boundary activity, the weak BCC efficiency is only sufficient but not necessary. Consequently, since a BCC-inefficient unit is not necessarily guaranteed to be an interior activity, it may lie on the frontier of $T^{DEA}$. As we will show in the next section, the reverse of Proposition 2.3 holds true when the data are positive.

## *2.3 DEA and Pareto efficiency based Classifications*

Depending on whether a strongly BCC-efficient DMU is either an extreme or a non-extreme point of $T^{DEA}$, it is called either extreme or non-extreme (Charnes et al., 1991). Therefore, from the equivalence of the strong BCC and Pareto efficiencies, we denote the sets of extreme and non-extreme efficient DMUs by $E$ and $E'$, respectively. However, due to non-equivalence of the weak BCC and Pareto efficiencies, $WE_I$ and $NE_I$ are used to denote the sets of weakly—but not strongly—BCC-efficient and BCC-inefficient DMUs, respectively. Moreover, $WE_p$ and $NE_p$ stand for the sets of weakly—but not strongly—Pareto-efficient and Pareto-inefficient DMUs, respectively.

From Definitions 2.1 and 2.3, the DEA and Pareto efficiency based classifications of the observed DMUs are thus derived respectively as

$$O = E \cup E' \cup WE_I \cup NE_I, \qquad (6)$$

$$O = E \cup E' \cup WE_P \cup NE_P. \qquad (7)$$

As previously mentioned, the reverse of Proposition 2.3 may not be true for non-negative data and, hence, classifications (6) and (7) may be different. In order to obtain classifications (6) and (7), one can apply the following *radial directional super-efficiency (RDSE)* model of Ray (2008):

$$\begin{aligned}
\min \quad & \beta - \varepsilon \left( \mathbf{1}_m^T \mathbf{s}^- + \mathbf{1}_s^T \mathbf{s}^+ \right) \\
\text{subject to} \quad & \\
& -\mathbf{X}_o \boldsymbol{\mu} - \mathbf{s}^- \geqq -\mathbf{x}_o - \beta \mathbf{g}^-, \\
& \mathbf{Y}_o \boldsymbol{\mu} - \mathbf{s}^+ \geqq \mathbf{y}_o - \beta \mathbf{g}^+, \\
& \mathbf{1}_{n-1}^T \boldsymbol{\mu} = 1, \\
& \boldsymbol{\mu} \geqq \mathbf{0}_{n-1}, \ \mathbf{s}^- \geqq \mathbf{0}_m, \ \mathbf{s}^+ \geqq \mathbf{0}_s,
\end{aligned} \quad (8)$$

where $\mathbf{g} = \left( \mathbf{g}^-, -\mathbf{g}^+ \right) \in \mathbb{R}_+^m \times \left( -\mathbb{R}_+^s \right)$ is a pre-specified direction vector along which $\mathbf{x}_o$ and $\mathbf{y}_o$ are, respectively, expanded and contracted.

**Proposition 2.4** (Dulá and Hickman 1997) Let $\mathbf{g} = \left( \mathbf{x}_o, \mathbf{0}_s \right)$ and let $\left( \beta^*, \boldsymbol{\mu}^*, \mathbf{s}^{-*}, \mathbf{s}^{+*} \right)$ be an optimal solution to model (8). Then,

(i) $(\mathbf{x}_o, \mathbf{y}_o) \in E$ if and only if $\beta^* > 0$ or model (8) is infeasible.

(ii) $(\mathbf{x}_o, \mathbf{y}_o) \in E'$ if and only if $\beta^* = 0$ and $\begin{pmatrix} \mathbf{s}^{-*} \\ \mathbf{s}^{+*} \end{pmatrix} = \mathbf{0}_{m+s}$.

(iii) $(\mathbf{x}_o, \mathbf{y}_o) \in WE_I$ if and only if $\beta^* = 0$ and $\begin{pmatrix} \mathbf{s}^{-*} \\ \mathbf{s}^{+*} \end{pmatrix} \geq \mathbf{0}_{m+s}$.

(iv) $(\mathbf{x}_o, \mathbf{y}_o) \in NE_I$ if and only if $\beta^* < 0$.

**Proposition 2.5** (Jahanshahloo et al. 2005) Let $\mathbf{g} = \left( \mathbf{1}_m, -\mathbf{1}_s \right)$ and let $\left( \beta^*, \boldsymbol{\Lambda}^*, \mathbf{s}^{-*}, \mathbf{s}^{+*} \right)$ be an optimal solution to model (8). Then,

(i) $(\mathbf{x}_o, \mathbf{y}_o) \in E$ if and only if $\beta^* > 0$ or model (8) is infeasible.

(ii) $(\mathbf{x}_o, \mathbf{y}_o) \in E'$ if and only if $\beta^* = 0$ and $\begin{pmatrix} \mathbf{s}^{-*} \\ \mathbf{s}^{+*} \end{pmatrix} = \mathbf{0}_{m+s}$.

(iii) $(\mathbf{x}_o, \mathbf{y}_o) \in WE_P$ if and only if $\beta^* = 0$ and $\begin{pmatrix} \mathbf{s}^{-*} \\ \mathbf{s}^{+*} \end{pmatrix} \geq \mathbf{0}_{m+s}$.

(iv) $(\mathbf{x}_o, \mathbf{y}_o) \in NE_P$ if and only if $\beta^* < 0$.

## *2.4 Classification based on both DEA and Pareto efficiencies*

As pointed out earlier, the BCC-inefficient units may be located on either $\partial T^{DEA}$ or $\text{int}(T^{DEA})$ in the presence of zero values in input–output data. Based on this fact, we partition $NE_I$ into the two disjoint groups $NW_I := NE_I \cap WE_P$ and $NN_I := NE_I \cap NE_P$. This separation results in the following input-oriented classification of the observed DMUs:

$$O = E \cup E' \cup WE_I \cup NW_I \cup NN_I. \tag{9}$$

As an illustration, consider units E and H in Fig. 2. As can be seen, both units are BCC inefficient since both of their inputs can be decreased radially. Moreover, while there exists no activity that strongly dominates E, all the points on the triangular ABC (except A) strongly dominate H. Hence, by Definition 2.3, the BCC-inefficient units E and H are weakly Pareto efficient and Pareto inefficient, respectively. Therefore, for $E, H \in NE_I$, we have $E \in NW_I$ and $H \in NN_I$.

**Remark 2.2** In an analogous manner, the use of the output-oriented BCC model yields the following output-oriented classification of the observed DMUs:

$$O = E \cup E' \cup WE_O \cup NW_O \cup NN_O, \tag{10}$$

where $NW_O := NE_O \cap WE_P$ and $NN_O := NE_O \cap NE_P$; here, $WE_O$ and $NE_O$ represent the sets of output-oriented weakly BCC-efficient and BCC-inefficient DMUs, respectively.

In order to gain both classifications (9) and (10), the following steps are required to be executed:

**Step 1** Classify the observed DMUs as in (6).

**Step 2** Reclassify the unit in $WE_I \cup NE_I$ into two groups: $WE_O$ and $NE_O$.

**Step 3** Reclassify the unit in $NE_I \cup NE_O$ into two groups: $WE_P$ and $NE_P$.

The execution of the above steps requires the first and second stages of model (8) to be run $n + Card(WE_I \cup NE_I) + Card(NE_I \cup NE_O)$ and $Card(E' \cup WE_I)$ times, respectively. In the immediately following section, we propose a dominance-based method to perform the aforementioned task by solving of a unique single-stage model $n$ times.

## 3 A unified dominance-based model

As the first step in developing our classification method in input-oriented case, we propose the following definition to relate the Farrell efficiency to the concept of dominance.

Note that the established relationships hold in the output orientation as well and are, hence, omitted.

**Definition 3.1** Let $(\mathbf{x}', \mathbf{y}')$ and $(\mathbf{x}, \mathbf{y})$ be two input–output combinations in $T^{DEA}$. Then, with an input orientation,

- $(\mathbf{x}', \mathbf{y}')$ *dominates* $(\mathbf{x}, \mathbf{y})$ in FGL sense, i.e., $(\mathbf{x}', \mathbf{y}') \succeq_I (\mathbf{x}, \mathbf{y})$, if and only if $\mathbf{x} >^+ \mathbf{x}'$ and $\mathbf{y}' \geqq \mathbf{y}$.

In the following lemma, we demonstrate that the definition of Farrell efficiency is based on the FGL dominance.

**Lemma 3.1** $DMU_o$ is

- BCC inefficient if and only if $(\mathbf{x}, \mathbf{y}) \succeq_I (\mathbf{x}_o, \mathbf{y}_o)$ holds for some $(\mathbf{x}, \mathbf{y}) \in T^{DEA}$.
- weakly BCC efficient if and only if $(\mathbf{x}, \mathbf{y}) \succeq_I (\mathbf{x}_o, \mathbf{y}_o)$ holds for no $(\mathbf{x}, \mathbf{y}) \in T^{DEA}$.

See Appendix for the proof.

Note that for positive data, the FGL and (weakly) Pareto dominances and, consequently, the weak BCC and weak Pareto efficiencies are equivalent by Definition 2.3 and Lemma 3.1.

Let us now consider the Pareto cone $K = \left\{ \begin{pmatrix} -\mathbf{s}^- \\ \mathbf{s}^+ \end{pmatrix} \middle| \mathbf{s}^- \in \mathbb{R}_+^m, \mathbf{s}^+ \in \mathbb{R}_+^s \right\}$. We define $\Omega_o$ as the intersection of the translated cone $\begin{pmatrix} \mathbf{x}_o \\ \mathbf{y}_o \end{pmatrix} + K$ with $T_o^{DEA}$, i.e., $\Omega_o = \left[ \begin{pmatrix} \mathbf{x}_o \\ \mathbf{y}_o \end{pmatrix} + K \right] \cap T_o^{DEA}$.

Then, this set specifies the set of possible activities in $T_o^{DEA}$ that all dominate $DMU_o$ and can be expressed equivalently as

$$\Omega_o = \left\{ \begin{pmatrix} \mathbf{x}_o - \mathbf{s}^- \\ \mathbf{y}_o + \mathbf{s}^+ \end{pmatrix} \middle| \mathbf{x}_o \geq \mathbf{X}_o \boldsymbol{\mu} + \mathbf{s}^-, \mathbf{Y}_o \boldsymbol{\mu} - \mathbf{s}^+ \geq \mathbf{y}_o, \mathbf{1}_{n-1}^T \boldsymbol{\mu} = 1, \boldsymbol{\mu} \geq \mathbf{0}_{n-1}, \mathbf{s}^- \geq \mathbf{0}_m, \mathbf{s}^+ \geq \mathbf{0}_s \right\}. \quad (11)$$

**Proposition 3.1** Let $n^+(\cdot)$ stands for the number of positive components of a vector and let $\bar{\Omega}_o = \Omega_o \setminus \left\{ \begin{pmatrix} \mathbf{x}_o \\ \mathbf{y}_o \end{pmatrix} \right\}$. Then,

(i) $(\mathbf{x}_o, \mathbf{y}_o) \in E$ if and only if $\Omega_o = \emptyset$.

(ii) $(\mathbf{x}_o, \mathbf{y}_o) \in E'$ if and only if $\bar{\Omega}_o = \varnothing$.

(iii) $(\mathbf{x}_o, \mathbf{y}_o) \in WE_I$ if and only if $\bar{\Omega}_o \neq \varnothing$ and $n^+(\mathbf{s}^-) < n^+(\mathbf{x}_o)$ for all $\begin{pmatrix} \mathbf{x}_o - \mathbf{s}^- \\ \mathbf{y}_o + \mathbf{s}^+ \end{pmatrix} \in \bar{\Omega}_o$.

(iv) $(\mathbf{x}_o, \mathbf{y}_o) \in NE_I$ if and only if $\bar{\Omega}_o \neq \varnothing$ and $n^+(\mathbf{s}^-) = n^+(\mathbf{x}_o)$ for some $\begin{pmatrix} \mathbf{x}_o - \mathbf{s}^- \\ \mathbf{y}_o + \mathbf{s}^+ \end{pmatrix} \in \bar{\Omega}_o$.

(v) $(\mathbf{x}_o, \mathbf{y}_o) \in WE_O$ if and only if $\bar{\Omega}_o \neq \varnothing$ and $n^+(\mathbf{s}^+) < n^+(\mathbf{y}_o)$ for all $\begin{pmatrix} \mathbf{x}_o - \mathbf{s}^- \\ \mathbf{y}_o + \mathbf{s}^+ \end{pmatrix} \in \bar{\Omega}_o$.

(vi) $(\mathbf{x}_o, \mathbf{y}_o) \in NE_O$ if and only if $\bar{\Omega}_o \neq \varnothing$ and $n^+(\mathbf{s}^+) = n^+(\mathbf{y}_o)$ for some $\begin{pmatrix} \mathbf{x}_o - \mathbf{s}^- \\ \mathbf{y}_o + \mathbf{s}^+ \end{pmatrix} \in \bar{\Omega}_o$.

(vii) $(\mathbf{x}_o, \mathbf{y}_o) \in WE_P$ if and only if $\bar{\Omega}_o \neq \varnothing$ and $n^+(\mathbf{s}^-) + n^+(\mathbf{s}^+) < m + s$ for all $\begin{pmatrix} \mathbf{x}_o - \mathbf{s}^- \\ \mathbf{y}_o + \mathbf{s}^+ \end{pmatrix} \in \bar{\Omega}_o$.

(viii) $(\mathbf{x}_o, \mathbf{y}_o) \in NE_P$ if and only if $\bar{\Omega}_o \neq \varnothing$ and $n^+(\mathbf{s}^-) + n^+(\mathbf{s}^+) = m + s$ for some $\begin{pmatrix} \mathbf{x}_o - \mathbf{s}^- \\ \mathbf{y}_o + \mathbf{s}^+ \end{pmatrix} \in \bar{\Omega}_o$.

See Appendix for the proof.

We now turn to visualize Proposition 3.1 with the help of a simple technology structure characterized by one input and one output.

**Example 3.1** Consider eight DMUs labeled as A, B, C, D, E, F and H, which construct the one-input and one-output technology. The empirical technological structures of $T_B^{DEA}$, $T_D^{DEA}$, $T_H^{DEA}$, $T_F^{DEA}$ and $T_G^{DEA}$ together with the corresponding translated cones are all exhibited in Fig.1 (a)–(e). The intersection of each technological structure with its corresponding translated cone is colored in yellow.

As can be seen in Fig. 1 (a), the intersection of $K+B$ with $T_B^{DEA}$ is empty, i.e., $\Omega_B = \varnothing$. This indicates that unit B is extreme efficient.

Fig. 1 (b) shows that the intersection of $K+D$ with $T_D^{DEA}$ is non-empty and is equal to the observed input–output vector itself, i.e., $\begin{pmatrix} x_D \\ y_D \end{pmatrix}$. This implies that unit D is non-extreme efficient.

From Fig. 1 (c), it can be seen that the intersection of $K+H$ with $T_H^{DEA}$ is the line segment joining A and H, and so $\bar{\Omega}_H \neq \varnothing$. Moreover, for all the points in $\bar{\Omega}_H$, it holds that $s^- = 0$ and $s^+ > 0$. This indicates that unit H is dominated but not strongly dominated. Therefore, H is weakly BCC efficient in input orientation, BCC inefficient in output orientation and weakly Pareto efficient.

Similarly, as depicted In Fig. 1 (d), the intersection of $K+F$ with $T_F^{DEA}$ is the line segment joining E and F, and thus $\bar{\Omega}_F \neq \varnothing$. In addition, $s^- > 0$ and $s^+ = 0$ hold for all the points in $\bar{\Omega}_F$. Hence, unit F is dominated but not strongly dominated. Therefore, F is BCC inefficient in input orientation, weakly BCC efficient in output orientation and weakly Pareto efficient.

It can be seen from Fig. 1 (e) that $\bar{\Omega}_G \neq \varnothing$. Moreover, $s^- > 0$ and $s^+ > 0$ for those points in $\Omega_G$ which are not located on the boundary of $K+G$. Thus, unit G is Pareto-inefficient.

In summary, we establish that the emptiness of $\Omega_o$ is both necessary and sufficient for $DMU_o$ to be extreme efficient. In the case of non-emptiness, we can also find out, with a little care, that the classification of $DMU_o$ relates primarily to the maximum values of $n^+(\mathbf{s}^-)$ and $n^+(\mathbf{s}^+)$ where $\begin{pmatrix} \mathbf{x}_o - \mathbf{s}^- \\ \mathbf{y}_o + \mathbf{s}^+ \end{pmatrix} \in \bar{\Omega}_o$. Based on these findings, we proceed to propose an LP model that helps us to determine whether the set $\Omega_o$ is empty; and if the answer is no, it then enables us to obtain an element of $\bar{\Omega}_o$ for which $n^+(\mathbf{s}^-)$ and $n^+(\mathbf{s}^+)$ are both maximum[4].

---

[4] Bertsimas and Tsitsiklis (1997) presented an LP problem to find out a feasible solution with the maximum number of positive components to a homogeneous system of equations (Exercise 3.27, pp. 136).

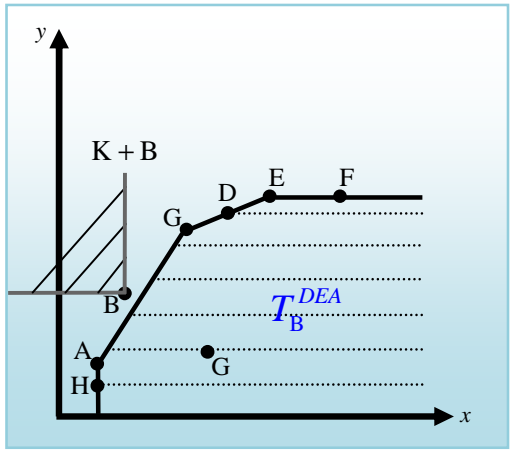
(a) extreme efficient

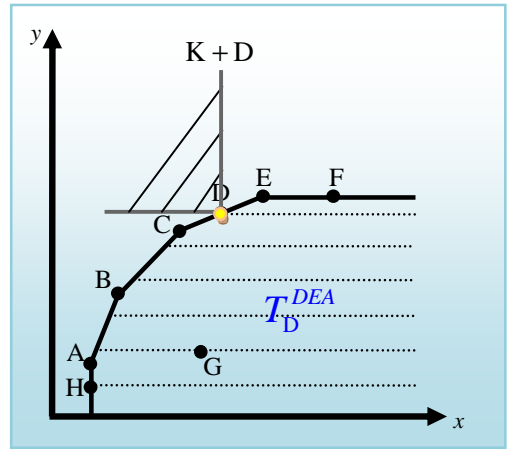
(b) non-extreme efficient

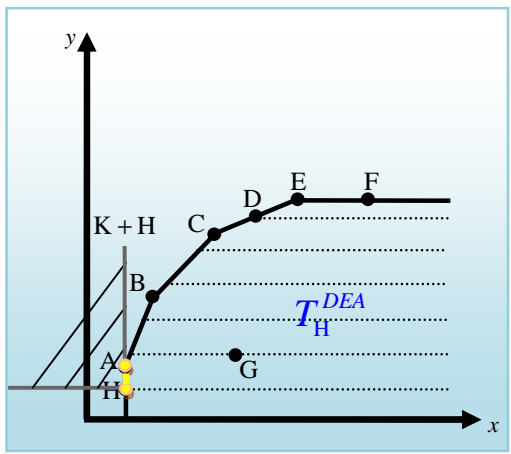
(c) weakly Pareto efficient,
input-oriented weakly BCC efficient,
output-oriented BCC inefficient

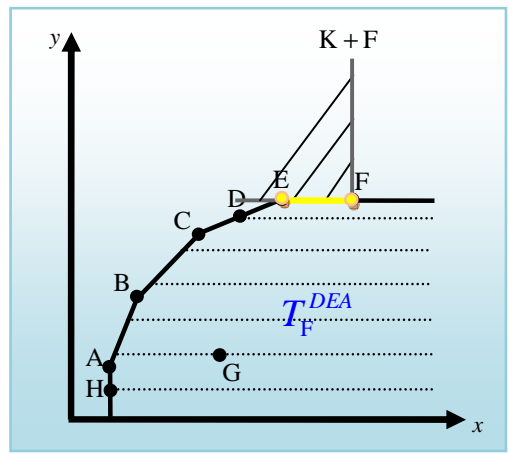
(d) weakly Pareto efficient,
input-oriented BCC inefficient,
output-oriented weakly BCC efficient

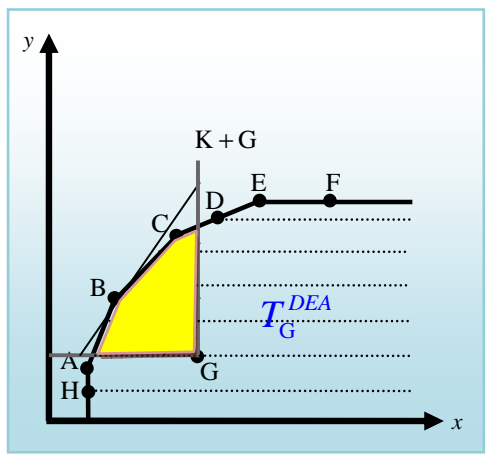
(e) Pareto inefficient

**Fig. 1** Technology structures with translated cones for selected DMUs

First, let $\begin{pmatrix} \mathbf{x}_o - \mathbf{s}^{-\prime} \\ \mathbf{y}_o + \mathbf{s}^{+\prime} \end{pmatrix} \in \Omega_o$. Then, (11) implies that $\mathbf{s}^{-\prime} \geq \mathbf{0}_m$, $\mathbf{s}^{+\prime} \geq \mathbf{0}_s$ and there exists $\boldsymbol{\mu}' \geq \mathbf{0}_{n-1}$ such that

$$\mathbf{x}_o \geq \mathbf{X}_o \boldsymbol{\mu}' + \mathbf{s}^{-\prime},\ \mathbf{Y}_o \boldsymbol{\mu}' - \mathbf{s}^{+\prime} \geq \mathbf{y}_o,\ \mathbf{1}_{n-1}^T \boldsymbol{\mu}' = 1. \tag{12}$$

We define

$$\sigma' := \max\left\{1, \min\left\{\frac{1}{s_i^{-\prime}}\middle|\, s_i^{-\prime} > 0\right\}, \min\left\{\frac{1}{s_r^{+\prime}}\middle|\, s_r^{+\prime} > 0\right\}\right\},\ \boldsymbol{\delta}' := \sigma' \boldsymbol{\mu}',$$

$$t_i^{-\prime} := \begin{cases} 1 & s_i^{-\prime} > 0, \\ 0 & s_i^{-\prime} = 0, \end{cases} \text{ for } i = 1,\ldots,m, \tag{13}$$

$$t_r^{+\prime} := \begin{cases} 1 & s_r^{+\prime} > 0, \\ 0 & s_r^{+\prime} = 0, \end{cases} \text{ for } r = 1,\ldots,s.$$

Then, (12) follows that $(\boldsymbol{\delta}', \mathbf{t}^{-\prime}, \mathbf{t}^{+\prime})$ is a feasible solution to the following system of constraints:

$$\begin{aligned} -\mathbf{X}_o \boldsymbol{\delta} + \left(\mathbf{1}_{n-1}^T \boldsymbol{\delta}\right) \mathbf{x}_o &\geq \mathbf{t}^-, \\ \mathbf{Y}_o \boldsymbol{\delta} - \left(\mathbf{1}_{n-1}^T \boldsymbol{\delta}\right) \mathbf{y}_o &\geq \mathbf{t}^+, \\ \mathbf{1}_{n-1}^T \boldsymbol{\delta} &\geq 1, \\ \boldsymbol{\delta} \geq \mathbf{0}_{n-1},\ \mathbf{1}_m \geq \mathbf{t}^- &\geq \mathbf{0}_m,\ \mathbf{1}_s \geq \mathbf{t}^+ \geq \mathbf{0}_s. \end{aligned} \tag{14}$$

As a direct consequence of the foregoing result, the following theorem holds.

**Proposition 3.2** For any $\begin{pmatrix} \mathbf{x}_o - \mathbf{s}^- \\ \mathbf{y}_o + \mathbf{s}^+ \end{pmatrix} \in \Omega_o$, there exists a feasible solution $(\boldsymbol{\delta}, \mathbf{t}^-, \mathbf{t}^+)$, as defined in (13), for system (14) such that $\mathbf{1}_m^T \mathbf{t}^- = n^+(\mathbf{s}^-)$ and $\mathbf{1}_s^T \mathbf{t}^+ = n^+(\mathbf{s}^+)$.

As shown above, the non-emptiness of $\Omega_o$ implies the feasibility of (14). To prove the converse, let $(\boldsymbol{\delta}', \mathbf{t}^{-\prime}, \mathbf{t}^{+\prime})$ be a feasible solution to (14). Then, the vector $\begin{pmatrix} \mathbf{x}_o - \mathbf{s}^{-\prime} \\ \mathbf{y}_o + \mathbf{s}^{+\prime} \end{pmatrix}$, defined by $\mathbf{s}^{-\prime} := \frac{1}{\mathbf{1}_{n-1}^T \boldsymbol{\delta}'} \mathbf{t}^{-\prime}$ and $\mathbf{s}^{+\prime} := \frac{1}{\mathbf{1}_{n-1}^T \boldsymbol{\delta}'} \mathbf{t}^{+\prime}$ with the corresponding intensity vector $\boldsymbol{\mu}' := \frac{1}{\mathbf{1}_{n-1}^T \boldsymbol{\delta}'} \boldsymbol{\delta}'$, is an element of $\Omega_o$, indicating that $\Omega_o \neq \emptyset$.

**Proposition 3.3** System (14) is feasible if and only if $\Omega_o \neq \emptyset$.

Proposition 3.3 states that the feasibility of system (14) is both necessary and sufficient for the non-emptiness of $\Omega_o$. Hence, by Proposition 3.1, the infeasibility of system (14) is both necessary and sufficient for DMU$_o$ to be extreme efficient, as stated in the following corollary.

**Corollary 3.1** $(\mathbf{x}_o, \mathbf{y}_o) \in E$ if and only if system (14) is infeasible.

With these results, we now propose the following LP model that maximizes $\mathbf{1}_m^T \mathbf{t}^- + \mathbf{1}_s^T \mathbf{t}^+$ subject to the constraints as stated in (14):

$$\begin{aligned}
\max \quad & \mathbf{1}_m^T \mathbf{t}^- + \mathbf{1}_s^T \mathbf{t}^+ \\
\text{subject to} \quad & \\
& -\mathbf{X}_o \boldsymbol{\delta} + \left(\mathbf{1}_{n-1}^T \boldsymbol{\delta}\right) \mathbf{x}_o \geqq \mathbf{t}^-, \\
& \mathbf{Y}_o \boldsymbol{\delta} - \left(\mathbf{1}_{n-1}^T \boldsymbol{\delta}\right) \mathbf{y}_o \geqq \mathbf{t}^+, \\
& \mathbf{1}_{n-1}^T \boldsymbol{\delta} \geq 1, \\
& \boldsymbol{\delta} \geqq \mathbf{0}_{n-1}, \; \mathbf{1}_m \geqq \mathbf{t}^- \geqq \mathbf{0}_m, \; \mathbf{1}_s \geqq \mathbf{t}^+ \geqq \mathbf{0}_s.
\end{aligned} \quad (15)$$

**Lemma 3.2** If model (15) is infeasible, then $(\mathbf{x}_o, \mathbf{y}_o) \in E$. Otherwise, both vectors $\mathbf{t}^{-*}$ and $\mathbf{t}^{+*}$ take on the values of 1 in all of their positive components.

See Appendix for the proof.

**Proposition 3.4** Let $(\boldsymbol{\delta}^*, \mathbf{t}^{-*}, \mathbf{t}^{+*})$ be an optimal solution to model (15) and define $\sigma^* := \mathbf{1}_{n-1}^T \boldsymbol{\delta}^*$. Then, $\mathbf{1}_m^T \mathbf{t}^{-*} = n^+(\mathbf{t}^{-*})$ and $\mathbf{1}_s^T \mathbf{t}^{+*} = n^+(\mathbf{t}^{+*})$. Furthermore, $\begin{pmatrix} \mathbf{x}_o - \dfrac{1}{\sigma^*} \mathbf{t}^{-*} \\ \mathbf{y}_o + \dfrac{1}{\sigma^*} \mathbf{t}^{+*} \end{pmatrix}$ is an element of $\Omega_o$ for which $n^+(\mathbf{t}^{-*}) + n^+(\mathbf{t}^{+*})$ is maximum.

See Appendix for the proof.

To this point, we have discriminated the extreme efficient DMUs. To classify the remaining DMUs under the feasibility of model (15), we now present the following corollary as a consequence of Proposition 3.4.

**Corollary 3.2** Let $\left(\boldsymbol{\delta}^{*}, \mathbf{t}^{-*}, \mathbf{t}^{+*}\right)$ be an optimal solution to model (15). Then,

(i) $\left(\mathbf{x}_{o}, \mathbf{y}_{o}\right) \in E'$ if and only if $\mathbf{1}_{m}^{T} \mathbf{t}^{-*} + \mathbf{1}_{s}^{T} \mathbf{t}^{+*} = 0$.

(ii) $\left(\mathbf{x}_{o}, \mathbf{y}_{o}\right) \in WE_{I}$ if and only if $0 < \mathbf{1}_{m}^{T} \mathbf{t}^{-*} + \mathbf{1}_{s}^{T} \mathbf{t}^{+*}$ and $\mathbf{1}_{m}^{T} \mathbf{t}^{-*} < n^{+}\left(\mathbf{x}_{o}\right)$.

(iii) $\left(\mathbf{x}_{o}, \mathbf{y}_{o}\right) \in NE_{I}$ if and only if $\mathbf{1}_{m}^{T} \mathbf{t}^{-*} = n^{+}\left(\mathbf{x}_{o}\right)$.

(iv) $\left(\mathbf{x}_{o}, \mathbf{y}_{o}\right) \in WE_{O}$ if and only if $0 < \mathbf{1}_{m}^{T} \mathbf{t}^{-*} + \mathbf{1}_{s}^{T} \mathbf{t}^{+*}$ and $\mathbf{1}_{s}^{T} \mathbf{t}^{+*} < n^{+}\left(\mathbf{y}_{o}\right)$.

(v) $\left(\mathbf{x}_{o}, \mathbf{y}_{o}\right) \in NE_{O}$ if and only if $\mathbf{1}_{s}^{T} \mathbf{t}^{+*} = n^{+}\left(\mathbf{y}_{o}\right)$.

(vi) $\left(\mathbf{x}_{o}, \mathbf{y}_{o}\right) \in WE_{P}$ if and only if $0 < \mathbf{1}_{m}^{T} \mathbf{t}^{-*} + \mathbf{1}_{s}^{T} \mathbf{t}^{+*} < m + s$.

(vii) $\left(\mathbf{x}_{o}, \mathbf{y}_{o}\right) \in NE_{P}$ if and only if $\mathbf{1}_{m}^{T} \mathbf{t}^{-*} + \mathbf{1}_{s}^{T} \mathbf{t}^{+*} = m + s$.

In summary, Corollaries 3.1 and 3.2 show that classifications (9) and (10) can be obtained by solving the single-stage model (15) $n$ times, one for each DMU.

From parts (ii), (iv) and (vi) of Corollary 3.2, the reverse of Proposition 2.3 holds true for positive data.

**Remark 3.1** In the presence of negative data[5], the Pareto efficiency based classification (7) can be made in three ways. The first one is to make the original data positive by adding some non-negative vector and then, to use the two-stage model (8) with $\mathbf{g} = \left(\mathbf{x}_{o}, \mathbf{0}_{s}\right)$ to the translated data. The second way is to apply this model directly to the original data by assigning $\mathbf{g} = \left(\mathbf{1}_{m}, -\mathbf{1}_{s}\right)$. The third way is to employ our single-stage model (15) directly to the original data. We have elaborated this important point through Example 4.2.

## 4 Numerical example

**Example 4.1** Consider a two-inputs and one-output technology characterized by eight hypothetical DMUs labeled as A–H. The observed input–output data are all exhibited in Table 1 and the resulting technology set is depicted in Fig. 2.

In order to illustrate the application of our proposed approach, we have solved model (15) for all the DMUs that are to be classified and have presented the results in Table 2. As shown in the second column of Table 2, model (15) is infeasible for units A, B and C, but feasible for the remaining ones. Hence, as can be seen in Fig. 2, units A, B and C are identified as extreme efficient by Corollary 3.1. Moreover, the remaining DMUs, for which model (15) is feasible, are classified based on both Corollary 3.2 and the optimal values $t_1^{-*}$, $t_2^{-*}$ and $t_1^{+*}$, which are all reported in the last three columns of Table 2. We observe that these values are all equal to one, as demonstrated in the proof of Lemma 3.2.

---

[5] While dealing with the estimation of a piecewise log-linear technology, one may encounter negative data since the log transformation of values less than 1 are always negative (Zarepisheh et al. 2010; Mehdiloozad et al. 2014). One may also refer to, e.g., Pastor and Ruiz (2007), Sahoo and Tone (2009) and Sahoo et al. (2012), among others, for several examples of applications with negative data.

**Table 1** Input-output data for eight DMUs

| DMU | $x_1$ | $x_2$ | $y$ |
|---|---|---|---|
| A | 0 | 1 | 1 |
| B | 2 | 1 | 2 |
| C | 0 | 2 | 2 |
| D | 0 | 1.5 | 1.5 |
| E | 0 | 4 | 1 |
| F | 2 | 1 | 1 |
| G | 4 | 4 | 2 |
| H | 4 | 4 | 1 |

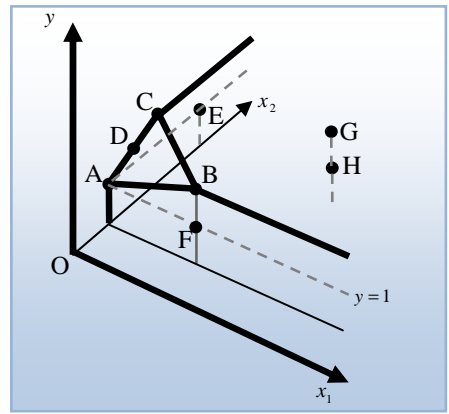

**Fig. 2** The technology set spanned by units A–H

**Table 2** The results for example

| DMU | Feasibility | $t_1^{-*}$ | $t_2^{-*}$ | $t_1^{+*}$ |
|---|---|---|---|---|
| A | I | -- | -- | -- |
| B | I | -- | -- | -- |
| C | I | -- | -- | -- |
| D | F | 0 | 0 | 0 |
| E | F | 0 | 1 | 1 |
| F | F | 1 | 0 | 1 |
| G | F | 1 | 1 | 0 |
| H | F | 1 | 1 | 1 |

In summary, complete classifications of the DMUs in both input and output orientations are all presented in Table 3.

**Table 3** The simultaneous classifications of DMUs in both input and output orientations

| DMU | Input-oriented | | | | | Output-oriented | | | | |
|---|---|---|---|---|---|---|---|---|---|---|
| | $E$ | $E'$ | $WE_I$ | $NW_I$ | $NN_I$ | $E$ | $E'$ | $WE_O$ | $NW_O$ | $NN_O$ |
| A | ✓ | | | | | ✓ | | | | |
| B | ✓ | | | | | ✓ | | | | |
| C | ✓ | | | | | ✓ | | | | |
| D | | ✓ | | | | | ✓ | | | |
| E | | | | ✓ | | | | | ✓ | |
| F | | | ✓ | | | | | | ✓ | |
| G | | | | ✓ | | | | ✓ | | |
| H | | | | | ✓ | | | | | ✓ |

**Example 4.2** Now, let us consider an example of a simple one-input and one-output technology that deals with both positive and negative data as described in Table 4. Via this example, we have shown how the single-stage model (15) can be applied for the Pareto-efficiency based classification of DMUs in the presence of negative data.

**Table 4** Input–output data for eight DMUs

| | A | B | C | D | E | F | G | H |
|---|---|---|---|---|---|---|---|---|
| Input | -6 | -6 | -5 | -4 | -2 | 1 | -4 | 2 |
| Output | -3 | 0 | 2 | 3 | 5 | 5 | 0 | 1 |

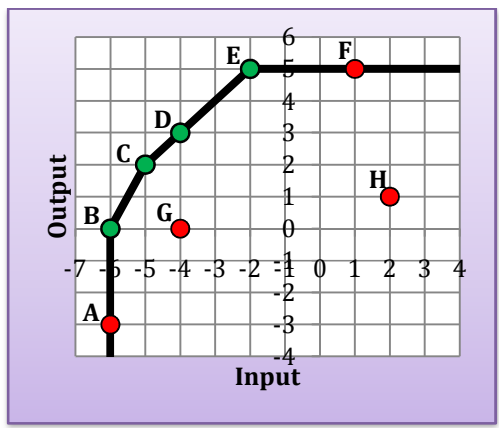

**Fig. 3** The production frontier

Fig. 3 displays the frontier spanned in the two-dimensional input–output space. From this figure, we observe that units A–F are boundary activities, amongst which units B, C and E are extreme efficient, unit D is non-extreme efficient, and units A and F are weakly Pareto efficient. Moreover, it can also be seen that units G and H are interior activities and are, hence, Pareto inefficient.

Table 5 presents the results obtained from model (15). As can be seen from the second column of Table 5, model (15) is infeasible for units B, C and E, thus confirming that these units are extreme efficient. For the remaining units, however, model (15) is feasible and the corresponding optimal values $t_1^{-*}$ and $t_1^{+*}$ are all reported in the last two columns.

**Table 5** The results for example

| DMU | Feasibility | $t_1^{-*}$ | $t_2^{-*}$ |
|---|---|---|---|
| A | F | 0 | 1 |
| B | I | -- | -- |
| C | I | -- | -- |
| D | F | 0 | 0 |
| E | I | -- | -- |
| F | F | 1 | 0 |
| G | F | 1 | 1 |
| H | F | 1 | 1 |

Now, applying Corollary 3.2 to the optimal values $t_1^{-*}$ and $t_1^{+*}$ yields the Pareto-efficiency based classification of eight DMUs, which are reported in Table 6.

**Table 6** The Pareto-efficiency based classification of DMUs

| DMU | E | E' | $WE_P$ | $NW_P$ |
|---|---|---|---|---|
| A |   |   | ✓ |   |
| B | ✓ |   |   |   |
| C | ✓ |   |   |   |
| D |   | ✓ |   |   |
| E | ✓ |   |   |   |
| F |   |   | ✓ |   |
| G |   |   |   | ✓ |
| H |   |   |   | ✓ |

## 5 Concluding remarks

In general, the classification of observed DMUs can be made with regard to either DEA efficiency or Pareto efficiency or both. While the DEA- and Pareto-efficiency based classifications have previously been addressed independently in the literature, no research has been published to date pertaining to the classification in terms of both DEA and Pareto efficiencies. The focus of the current study is, therefore, to accomplish this task in an efficient manner using an alternative dominance-based method.

It is well known that both DEA and Pareto efficiencies are either weak or strong. Though the strong DEA and Pareto efficiencies are equivalent, the statement is generally not true for weak efficiency. Precisely, the weak DEA efficiency is sufficient—but not necessary—for achieving the Pareto efficiency. Since the boundary activities are characterized by the concept of weak Pareto efficiency, a DEA-inefficient DMU can therefore be either an interior or a boundary activity. Based on these facts, we reclassify the observed DMUs into five groups: extreme efficient, non-extreme efficient, weakly—but not strongly—DEA efficient, boundary DEA inefficient and interior DEA inefficient.

To make our proposed classification, we employed the concept of FGL dominance to show that, like the Pareto efficiency, the DEA efficiency is also a dominance-based concept. We then propose a single-stage linear programming model that enables us to simultaneously classify the observed DMUs in both input and output orientations. Our model is also able to make the Pareto-efficiency based classification in the presence of negative data. Moreover, by identifying all possible improvements in the inputs and outputs, the proposed model enables a decision maker to deal with his/her preference issue of which specific input to reduce or which specific output to increase for a weakly Pareto-efficient DMU in order to improve its overall performance. Based on this property, we point to avenues for future research as to how our proposed model can be implemented for target setting and sensitivity analysis.

## Appendix

**Proof of Lemma 3.1**

Since part (ii) is an immediate consequence of part (i), we prove only part (i).

Let $(\mathbf{x}_o, \mathbf{y}_o) \in NE_I$. By Definition 2.1, there exists $\hat{\theta} < 1$ such that $(\hat{\theta}\mathbf{x}_o, \mathbf{y}_o) \in T^{DEA}$. Since $\hat{\theta} < 1$, we have $\mathbf{x}_o >^+ \hat{\theta}\mathbf{x}_o$. Hence, by Definition 3.1, it follows that $(\hat{\theta}\mathbf{x}_o, \mathbf{y}_o) \succeq_I (\mathbf{x}_o, \mathbf{y}_o)$.

Conversely, let $(\hat{\mathbf{x}}, \hat{\mathbf{y}}) \in T^{DEA}$ dominates $(\mathbf{x}_o, \mathbf{y}_o)$ in the FGL sense. By Definition 3.1, we have $\mathbf{x}_o >^+ \hat{\mathbf{x}}$, $\hat{\mathbf{y}} \geq \mathbf{y}$. We define $\hat{\theta} := \max\left\{\frac{\hat{x}_i}{x_{io}} \mid x_{io} > 0\right\}$. Then, $\hat{\theta} < 1$ and $\hat{\theta}\mathbf{x}_o \geq \hat{\mathbf{x}}$.

Since $(\hat{\mathbf{x}},\hat{\mathbf{y}}) \in T^{DEA}$, as per the free disposability assumption, $(\hat{\theta}\mathbf{x}_o,\mathbf{y}_o) \in T^{DEA}$. Hence, the optimal value of $\theta^*$ in model (5) is less than one, indicating $(\mathbf{x}_o,\mathbf{y}_o) \in NE_I$. ∎

**Proof of Proposition 3.1**

Note that parts (iii), (v) and (vii) are immediate consequences of parts (iv), (vi) and (viii), respectively. Moreover, the proof of part (vi) is similar to that of part (iv). Therefore, we are left with in proving only parts (i), (ii), (iv) and (viii).

**Part (i)** Let $(\mathbf{x}_o,\mathbf{y}_o) \in E$ and assume by contradiction that $\begin{pmatrix} \mathbf{x}_o - \hat{\mathbf{s}}^- \\ \mathbf{y}_o + \hat{\mathbf{s}}^+ \end{pmatrix} \in \Omega_o$. Let $\hat{\boldsymbol{\mu}}$ be the intensity variable corresponding to this element in (11). For any direction vector $\mathbf{g}$, $(\beta',\boldsymbol{\mu}',\mathbf{s}^{-\prime},\mathbf{s}^{+\prime}) := (0,\hat{\boldsymbol{\mu}},\hat{\mathbf{s}}^-,\hat{\mathbf{s}}^+)$ is then a feasible solution to model (8). Since the first phase of model (8) is a minimization problem, the optimal value of $\beta$ is non-positive, which is a contradiction by Proposition 2.4.

In a similar way, the converse can be proved by the way of contradiction and with the help of Proposition 2.4.

**Part (ii)** Let $(\mathbf{x}_o,\mathbf{y}_o) \in E'$. Then, part (ii) of Proposition 2.4 follows that $\left\{\begin{pmatrix} \mathbf{x}_o \\ \mathbf{y}_o \end{pmatrix}\right\} \subseteq \Omega_o$. To prove the equality, assume by contradiction that $\begin{pmatrix} \mathbf{x}_o - \hat{\mathbf{s}}^- \\ \mathbf{y}_o + \hat{\mathbf{s}}^+ \end{pmatrix} \in \Omega_o$ such that $\begin{pmatrix} \hat{\mathbf{s}}^- \\ \hat{\mathbf{s}}^+ \end{pmatrix} \geq \mathbf{0}_{m+s}$. Let $\hat{\boldsymbol{\mu}}$ be the intensity variable corresponding to this element in (11). For any direction vector $\mathbf{g}$, $(\beta',\boldsymbol{\mu}',\mathbf{s}^{-\prime},\mathbf{s}^{+\prime}) := (0,\hat{\boldsymbol{\mu}},\hat{\mathbf{s}}^-,\hat{\mathbf{s}}^+)$ is then a feasible solution to model (8). Since the second phase of model (8) is a maximization problem, the optimal sum of slacks is positive, which is a contradiction by Proposition 2.4.

In a similar way, the converse can be proved by the way of contradiction and with the help of Proposition 2.4.

**Part (iv)** By Lemma 3.1, $(\mathbf{x}_o,\mathbf{y}_o) \in NE_I$ if and only if there exists $(\hat{\mathbf{x}},\hat{\mathbf{y}}) \in T^{DEA}$ such that $\mathbf{x}_o >^+ \hat{\mathbf{x}}$ and $\mathbf{y}' \geq \mathbf{y}_o$. We define $\hat{\mathbf{s}}^- := \mathbf{x}_o - \hat{\mathbf{x}}$ and $\hat{\mathbf{s}}^+ := \hat{\mathbf{y}} - \mathbf{y}_o$. Then, $(\mathbf{x}_o,\mathbf{y}_o) \in NE_I$ if and only if (i) $\hat{\mathbf{s}}^- \geq \mathbf{0}_m$ and $\hat{\mathbf{s}}^+ \geq \mathbf{0}_s$, and (ii) $\hat{s}_i^- > 0 \Leftrightarrow x_{io} > 0$; or, if and only if $\begin{pmatrix} \mathbf{x}_o - \hat{\mathbf{s}}^- \\ \mathbf{y}_o + \hat{\mathbf{s}}^+ \end{pmatrix} \in \bar{\Omega}_o$ and $n^+(\hat{\mathbf{s}}^-) = n^+(\mathbf{x}_o)$.

**Part (viii)** By Definition 2.3, $(\mathbf{x}_o, \mathbf{y}_o) \in NE_P$ if and only if there exists $(\hat{\mathbf{x}}, \hat{\mathbf{y}}) \in T^{DEA}$ that strongly dominates $(\mathbf{x}_o, \mathbf{y}_o)$. We define $\hat{\mathbf{s}}^- := \mathbf{x}_o - \hat{\mathbf{x}}$ and $\hat{\mathbf{s}}^+ := \hat{\mathbf{y}} - \mathbf{y}_o$. Then, $(\mathbf{x}_o, \mathbf{y}_o) \in NE_P$ if and only if $\begin{pmatrix} \hat{\mathbf{s}}^- \\ \hat{\mathbf{s}}^+ \end{pmatrix} > \mathbf{0}_{m+s}$; or, if and only if $\begin{pmatrix} \mathbf{x}_o - \hat{\mathbf{s}}^- \\ \mathbf{y}_o + \hat{\mathbf{s}}^+ \end{pmatrix} \in \bar{\Omega}_o$ and $n^+(\hat{\mathbf{s}}^-) + n^+(\hat{\mathbf{s}}^+) = m + s$. ∎

**Proof of Lemma 3.2**

First, assume that model (15) is infeasible. Then, system (14) is infeasible and, by Corollary 3.1, $(\mathbf{x}_o, \mathbf{y}_o)$ is an extreme-efficient DMU.

Now, let model (15) be feasible. Then, this model has an optimal solution since its objective function is upper bounded by $m+s$. Hence, let $(\boldsymbol{\delta}^*, \mathbf{t}^{-*}, \mathbf{t}^{+*})$ be an optimal solution to model (15). We claim that the positive components of the vectors $\mathbf{t}^{-*}$ and $\mathbf{t}^{+*}$ are all equal to one. Since the proofs are similar for these vectors, we prove the assertion only for $\mathbf{t}^{-*}$.

By way of contradiction, assume that $0 < t_g^{-*} < 1$ for some $g \in \{1,...,m\}$. Dividing both sides of the constraints of (15) at optimality by $t_g^{-*}$ yields

$$\mathbf{X}_o \left( \frac{1}{t_g^{-*}} \boldsymbol{\delta}^* \right) + \frac{1}{t_g^{-*}} \mathbf{t}^{-*} \leq \left( \frac{\mathbf{1}_{n-1}^T \boldsymbol{\delta}^*}{t_g^{-*}} \right) \mathbf{x}_o,$$

$$\mathbf{Y}_o \left( \frac{1}{t_g^{-*}} \boldsymbol{\delta}^* \right) - \frac{1}{t_g^{-*}} \mathbf{t}^{+*} \geq \left( \frac{\mathbf{1}_{n-1}^T \boldsymbol{\delta}^*}{t_g^{-*}} \right) \mathbf{y}_o, \qquad (16)$$

$$\frac{\mathbf{1}_{n-1}^T \boldsymbol{\delta}^*}{t_g^{-*}} \geq 1,$$

Then, according to (16), the vector $(\boldsymbol{\delta}', \mathbf{t}^{-\prime}, \mathbf{t}^{+\prime})$ defined by

$$\boldsymbol{\delta}' := \frac{1}{t_g^{-*}} \boldsymbol{\delta}^*,$$

$$t_i^{-\prime} := \min\left\{1, \frac{t_i^{-*}}{t_g^{-*}}\right\}, \text{ for } i = 1,...,m, \qquad (17)$$

$$t_r^{+\prime} := \min\left\{1, \frac{t_r^{+*}}{t_g^{-*}}\right\}, \text{ for } r = 1,...,s,$$

is a feasible solution to model (15) whose objective function value is strictly greater than $\mathbf{1}_m^T \mathbf{t}^{-*} + \mathbf{1}_s^T \mathbf{t}^{+*}$. This contradicts the optimality of $(\boldsymbol{\delta}^*, \mathbf{t}^{-*}, \mathbf{t}^{+*})$ and proves our claim. ∎

**Proof of Proposition 3.4**

According Lemma 3.2, it is obvious that $\mathbf{1}_m^T \mathbf{t}^{-*} = n^+(\mathbf{t}^{-*})$ and $\mathbf{1}_s^T \mathbf{t}^{+*} = n^+(\mathbf{t}^{+*})$. Moreover, as in the paragraph below proposition 3.2, it is straightforward to show that

$$\begin{pmatrix} \mathbf{x}_o - \dfrac{1}{\sigma^*}\mathbf{t}^{-*} \\ \mathbf{y}_o + \dfrac{1}{\sigma^*}\mathbf{t}^{+*} \end{pmatrix} \in \Omega_o .$$

Let $\begin{pmatrix} \mathbf{x}_o - \hat{\mathbf{s}}^- \\ \mathbf{y}_o + \hat{\mathbf{s}}^+ \end{pmatrix} \in \Omega_o$ for which $n^+(\hat{\mathbf{s}}^-) + n^+(\hat{\mathbf{s}}^+)$ is maximum. Then, we have $\mathbf{1}_m^T \mathbf{t}^{-*} + \mathbf{1}_s^T \mathbf{t}^{+*} = n^+(\mathbf{t}^{-*}) + n^+(\mathbf{t}^{+*}) \leq n^+(\hat{\mathbf{s}}^-) + n^+(\hat{\mathbf{s}}^+)$. From Proposition 3.2, we also know that there exists a feasible solution $(\hat{\boldsymbol{\delta}}, \hat{\mathbf{t}}^-, \hat{\mathbf{t}}^+)$ for model (15) such that $\mathbf{1}_m^T \hat{\mathbf{t}}^- + \mathbf{1}_s^T \hat{\mathbf{t}}^+ = n^+(\hat{\mathbf{s}}^-) + n^+(\hat{\mathbf{s}}^+)$. Therefore, the proof is complete by the fact that model (15) is a maximization LP problem. ∎